\documentclass[11pt]{amsart}
\usepackage[a4paper,margin=1in]{geometry}
\usepackage[dvipsnames]{xcolor}
\usepackage{parskip}
\usepackage{microtype}
\usepackage{enumitem}
\usepackage{tikz}
\usetikzlibrary{calc,arrows.meta}

\tikzset{
 blueedge/.style={draw=RoyalBlue,line width=1pt},
 rededge/.style={draw=BrickRed,dashed,line width=1pt},
 greenedge/.style={draw=ForestGreen,densely dotted,
                   line width=1.35pt},
 allowededge/.style={draw=black!70,line width=.5pt,
                     double=white,double distance=.9pt},
 vertex/.style={circle,fill=black,inner sep=1.4pt},
 proofarrow/.style={-{Latex[length=2mm]},
                    draw=black!65,line width=.7pt}
}

\usepackage[
 colorlinks=true,
 linkcolor=MidnightBlue,
 urlcolor=MidnightBlue,
 citecolor=MidnightBlue
]{hyperref}
\newcommand{\GRregular}[4]{%
  \pgfmathtruncatemacro{\GRlast}{#2-1}%
  \foreach \GRi in {0,...,\GRlast}
    \coordinate (#1\GRi) at ({#4-360*\GRi/#2}:#3);%
}

\newtheorem{theorem}{Theorem}
\newtheorem{lemma}[theorem]{Lemma}
\newtheorem{corollary}[theorem]{Corollary}

\newcommand{\gr}{\operatorname{gr}_3^2}
\newcommand{\g}[3]{\gr(K_{#1},K_{#2},K_{#3})}
\newcommand{\BR}{\mathrm{BR}}
\newcommand{\BG}{\mathrm{BG}}
\newcommand{\RG}{\mathrm{RG}}
\newcommand{\wbr}{\omega_{\BR}}
\newcommand{\wbg}{\omega_{\BG}}
\newcommand{\wrg}{\omega_{\RG}}
\newcommand{\case}[1]{\smallskip\noindent\textbf{#1.}}

\title[The weakened Gallai--Ramsey number of $K_7$]{The weakened Gallai--Ramsey number $\operatorname{gr}_3^2(K_7)=17$}

\author[Aggarwal et al.]{Mudit Aggarwal$^{*,\dagger}$}
\author[]{Andrew Alexander$^{*}$}
\author[]{Kyle Bierly$^{*}$}
\author[]{Jeffrey Dawson$^{*}$}
\thanks{$^{*}$Department of Mathematics, The University of British Columbia,
1984 Mathematics Road, Vancouver, BC V6T 1Z2, Canada.}
\thanks{$^{\dagger}$TU Wien, Karlsplatz 13, 1040 Vienna, Austria.}
\thanks{Corresponding author: Mudit Aggarwal, \texttt{mudit.aggarwal@tuwien.ac.at}.}


\begin{document}

\begin{abstract}
We study when Gallai three-colorings of complete graphs (edge-colorings with no rainbow triangle) must contain a complete subgraph whose edges use at most two colors. We prove that the weakened Gallai--Ramsey number $\operatorname{gr}_3^2(K_7)$ equals $17$: every Gallai three-coloring of $K_{17}$ contains a two-colored $K_7$, and we construct a Gallai three-coloring of $K_{16}$ that does not. In the course of the proof, we determine the previously unknown off-diagonal values $g(4,6,6)=9$, $g(5,5,7)=11$, and $g(6,6,7)=13$, where $g(a,b,c)$ denotes the least $N$ such that every Gallai three-coloring of $K_N$ contains a blue-red $K_a$, a blue-green $K_b$, or a red-green $K_c$.
\end{abstract}

\maketitle

\section{Introduction}
A $t$-coloring of a graph is a map (not necessarily surjective) from the edge set of the graph to a set of $t$ colors. A $t$-coloring of a graph is a \emph{Gallai coloring} if all triangles in the graph use at most two colors. It is also commonly said that such a coloring \emph{avoids rainbow triangles}.

Denoted by $\operatorname{gr}_{t}^s(G)$, the \emph{weakened Gallai--Ramsey number} is the smallest positive integer $N$ such that every Gallai $t$-coloring of the complete graph $K_N$ contains a subgraph isomorphic to $G$ whose edges use at most $s$ colors, where $1 \leq s \leq t$. In this paper, we only consider the case $s=2$ and $t=3$, with $G$ a complete graph. Throughout, \emph{two-colored} means using at most two colors; in particular, a monochromatic subgraph is two-colored.

Gallai colorings originate in Gallai's work~\cite{gallai}; their Ramsey-type study was initiated by Gy\'arf\'as and Simonyi~\cite{gs}. Ramsey numbers for subgraphs using at most $s$ of $t$ colors go back to Chung, Chung, and Liu and to Chung and Liu~\cite{ccl, cl}. Harborth and M\"oller~\cite{hm} studied weakened Ramsey numbers without the Gallai restriction. Beam and Budden~\cite{bb} introduced weakened Gallai--Ramsey numbers and proposed the off-diagonal version. In the three-color arguments, we consider the colors blue, red, and green. Since we only consider $s=2$ and $t=3$, we further write
\[
 g(a,b,c)=\g{a}{b}{c},
\]
which denotes the smallest positive integer $N$ such that every Gallai three-coloring of the complete graph $K_N$ contains at least one of the following:
\begin{enumerate}
    \item A blue-red $K_a$.
    \item A blue-green $K_b$.
    \item A red-green $K_c$.
\end{enumerate}
We call $g(a,b,c)$ \emph{off-diagonal} when $a$, $b$, and $c$ are not all equal. It immediately follows from the definition that $\gr(K_n) = g(n,n,n)$.

A foundational result of Gallai~\cite{gallai}, in the form given by Gy\'arf\'as and Simonyi~\cite[Theorem~A]{gs}, expresses every nontrivial Gallai coloring as a substitution into a two-colored complete graph. A formal statement follows.

\begin{theorem}[{\cite{gallai, gs}}]\label{thm:partition}
    For every nontrivial Gallai coloring of a complete graph, the vertex set of the graph can be partitioned into nonempty blocks $X_1,\ldots,X_r$, $r\ge2$. Every pair of blocks is joined by edges of one color, and at most two colors appear between the blocks. Each block is itself Gallai colored.
\end{theorem}
The decomposition defines the \emph{base graph} $H$: the two-colored complete graph on $\{1,\ldots,r\}$ in which the edge $ij$ has the color of the edges between $X_i$ and $X_j$. We say that $X_i$ and $X_j$ are \emph{joined} in this color, and that $X_j$ is a blue (red, green) \emph{neighbor} of $X_i$ if they are joined in blue (red, green). A block is \emph{all-blue}, or simply \emph{blue}, if every edge inside it is blue; in particular, a single vertex is an all-blue block (and similarly for red and green). Theorem~\ref{thm:partition} is used widely in proving bounds on values for graph classes (see, for example, \cite{jm, wz}).

Chung and Graham~\cite{cg} determined the following values for monochromatic triangles (a shorter proof appears in~\cite{gsss}):
\[
 \operatorname{gr}_t^1(K_3)=
 \begin{cases}5^{t/2}+1,&t\text{ even},\\
 2\cdot5^{(t-1)/2}+1,&t\text{ odd}.
 \end{cases}
\]

A result of Erd\H{o}s, Simonovits, and S\'os~\cite{ess} gives $\operatorname{gr}_t^s(K_{s+1})=s+1$ for $1\le s<t$ (see also \cite[Theorem~9]{jm}). In particular, $\operatorname{gr}_t^2(K_3)=3$ for $t\ge3$.

Beam and Budden~\cite[Theorem~6]{bb} proved that $\operatorname{gr}_t^{t-1}(K_{t+1})=t+2$ for $t\ge3$; in particular, $\gr(K_4)=5$. Budden and Wimbish proved that $\operatorname{gr}_t^2(K_4)=t+2$ for $t\ge3$ and that $\gr(K_5)=9$~\cite[Theorems~2.1--2.4]{bw}. They further proved that $\operatorname{gr}_t^2(K_5)\ge2^t+1$ for $t\ge3$, with equality for $t=3$~\cite[Theorems~2.2 and~2.3]{bw}. Their Conjecture~4.1, asserting equality for all $t\ge3$, was proved by Li, Broersma, and Wang~\cite[Theorem~1.11]{lbw}, who also established~\cite[Theorem~1.8]{lbw}
\[
 \operatorname{gr}_t^{m-2}(K_m)=t+2\quad(m\ge4,\ t\ge m-2).
\]

Using the Gallai analogue of the Erd\H{o}s--Gy\'arf\'as problem~\cite{eg}, it was shown~\cite{lbw} that
\[
 \operatorname{gr}_t^s(K_m)\le2^{\,2t(m-2)/s+1}
 \quad(m\ge3,\ 1\le s\le m-2,\ t\ge s).
\]
Fox, Grinshpun, and Pach~\cite[Theorem~1.1]{fgp} show that the largest two-colored clique guaranteed in every Gallai three-coloring of $K_N$ has order $\Theta(N^{1/3}(\log N)^2)$. Inverting their theorem yields
\[
 \gr(K_m)=\Theta \left( \frac{m^3}{(\log m)^6} \right) \quad (m\to\infty).
\]

Jakhar and Moun~\cite{jm} prove results for weakened Gallai--Ramsey numbers for several graph classes, not limited to complete graphs. Additionally, by constructing a Gallai three-coloring of $K_{11}$, they prove $\gr(K_6)\ge12$, and conjecture that equality holds.

Wang and Zhang~\cite[Theorem~1.3]{wz} resolve this conjecture by proving $\gr(K_6)=12$. In our notation, their auxiliary results~\cite[Lemmas~2.5--2.7]{wz} include
\[
 g(4,4,6)\le7,\qquad g(4,5,6)\le10,\qquad g(5,5,6)\le11.
\]
For more colors, they also establish
\[
 \operatorname{gr}_t^2(K_6)\ge
 \begin{cases}
  5^{t/2}+1,&t\text{ even},\\
  2\cdot5^{(t-1)/2}+t-1,&t\text{ odd},
 \end{cases}
\]
for every $t\ge2$, and conjecture equality~\cite[Conjecture~1.5]{wz}.

Budden and DeJonge~\cite{bd} develop the off-diagonal problem, determining $g(3,j,n)$ for $3\le j\le9$ and $n\ge j$, and the families $g(4,4,n)$ and $g(4,5,n)$. Their value $g(4,5,6)=8$ sharpens the bound $g(4,5,6)\le10$ used by Wang and Zhang. The specific cases that we require for our proofs are listed in Section~\ref{sec:asym}.

The aim of this paper is to prove that $\gr(K_7)=17$. To arrive at this result, we first collect known off-diagonal values and prove $g(4,6,6)=9$ (Theorem~\ref{thm:466}) and $g(5,5,7)=11$ (Theorem~\ref{thm:557}) in Section~\ref{sec:asym}. We then prove $g(6,6,7)=13$ (Theorem~\ref{thm:667}) in Section~\ref{sec:667}. Collectively, these yield Theorem~\ref{thm:k7} in Section~\ref{sec:k7}.

\section{Preliminaries}\label{sec:prelim}
Consider a complete graph \(X\) with a coloring on its edges, and a nonempty set of colors $S$. Let $\omega_S(X)$ be the largest order (number of vertices) of a clique in $X$ whose edge colors belong to $S$, and write $|X|$ for the order of $X$. If \(X\) is a single vertex, then \(\omega_S(X) = 1\) for every nonempty \(S\). We sometimes omit the argument $X$ when the graph is clear.

We now limit ourselves to three colors, and let $\BR=\{\text{blue},\text{red}\}$, $\BG=\{\text{blue},\text{green}\}$, and $\RG=\{\text{red},\text{green}\}$. For example, $\wbr(X)$ is the largest order of a clique in $X$ containing no green edge, and if $X$ is an all-blue $K_m$, then $\wbr(X)=\wbg(X)=m$ and $\wrg(X)=1$.

Given $a$, $b$, and $c$, we call a blue-red $K_a$, a blue-green $K_b$, and a red-green $K_c$ the \emph{targets} for $g(a,b,c)$. A coloring $X$ avoids the targets precisely when
\[
 \wbr(X)<a,\qquad \wbg(X)<b,\qquad \wrg(X)<c,
\]
and $g(a,b,c)>N$ precisely when some Gallai three-coloring of $K_N$ avoids them. This is the condition we work with for most of the paper. Permuting the colors shows that $g(a,b,c)$ is symmetric in $a$, $b$, and $c$; we use this freely, and remind the reader when we do.

We prove each equality $g(a, b, c) = n$ in two parts. For the upper bound $g(a,b,c)\le n$, we take a Gallai three-coloring of $K_n$ that avoids all three targets and derive a contradiction. For the lower bound $g(a,b,c)\ge n$, we construct an explicit Gallai three-coloring of $K_{n-1}$ that avoids all three targets.

For the upper bound, we typically make use of the decomposition given by Theorem~\ref{thm:partition}. We always choose a decomposition with the fewest possible blocks, and refer to such a decomposition as \emph{minimal}. Let $r$ be its number of blocks. If $r\ge3$, then the base graph has no \emph{monochromatic cut}, that is, no partition of its vertices into two nonempty parts with all edges between the parts of one color, since merging the blocks on each side would give a decomposition with two blocks. Every two-colored triangle has a vertex whose two edges share a color, so $r\ne3$. Hence $r=2$ or $r\ge4$ (see also \cite[Lemma~3.1]{mn}), and when $r\ge4$, every block has a neighbor in each base color.

We write $R(a,b)$ for the classical two-color Ramsey number. We require the classical value $R(3, 3)=6$ (see, for instance, \cite{gg}). In the proofs of results from \cite{bd} mentioned in Section~\ref{sec:asym}, the values $R(3,4)=9$ and $R(4, 4) = 18$ are also used.

For our results, we also use the following diagonal values for weakened Gallai--Ramsey numbers:
\[\gr(K_4)=5, \qquad \gr(K_5)=9, \qquad \gr(K_6)=12.\]
These are proved in \cite[Theorem~2.1]{bw}, \cite[Theorem~2.3]{bw}, and \cite[Theorem~1.3]{wz}, respectively.

The following observation is used in our proofs multiple times. The partition also appears in \cite[Lemma~2.2]{wz}; we give a short proof here for completeness.
\begin{lemma}\label{lem:triangle}
Let $X$ be a Gallai three-colored complete graph with no blue-red triangle. Then the following hold.
\begin{enumerate}[label=(\arabic*), ref=\thetheorem(\arabic*)]
\item\label{lem:tri1} No vertex is incident with both a blue and a red edge.
\item\label{lem:tri2} The vertices of $X$ can be partitioned into sets $A$ and $B$ such that every edge within $A$ is red or green, every edge within $B$ is blue or green, and every edge between $A$ and $B$ is green; either set may be empty.
\item\label{lem:tri3} If, in addition, $|X|\ge3$ and $X$ has no blue-green triangle, then $X$ has no blue edge, that is, $X$ is red-green.
\end{enumerate}
\end{lemma}
\begin{proof}
For (1), a blue and a red edge meeting at a vertex would complete to either a blue-red triangle or a rainbow triangle.

For (2), we put every vertex incident with a blue edge in $B$ and all remaining vertices in $A$. Vertices in $A$ are not incident with blue edges, and by (1), vertices in $B$ are not incident with red edges. Thus, all edges between $A$ and $B$ are green.

For (3), since $|X|\ge3$, every blue edge lies in a triangle, which, being Gallai, would be blue-red or blue-green, contrary to the hypotheses.
\end{proof}

We further note our stylistic choices in the figures that follow. For accessibility, we use \emph{solid blue}, \emph{dashed red}, and \emph{dotted green} edges. Edges that may take more than one color are drawn as \emph{hollow gray} double lines, with the possible colors given in the text or caption.

\section{Initial off-diagonal bounds} \label{sec:asym}

In this section, we prove some off-diagonal bounds that will be needed later on, and collect known results (primarily by Budden and DeJonge~\cite{bd}) in our notation.

\begin{theorem}[{\cite[Theorem~3.4]{bd}}]\label{thm:34n}
If $n\ge5$, then $g(3,4,n)=\min\{n,9\}$.
Also, $g(3,4,4)=5$.
\end{theorem}

We use the following consequence frequently.
\begin{corollary}\label{cor:both4}
Let $X$ be a Gallai three-colored complete graph with $|X| \geq 5$ that uses all three colors. If $\wbr(X)\le2$ (no blue-red triangle), then $\wbg(X),\wrg(X)\ge4$.
\end{corollary}
\begin{proof}
The coloring contains a blue and a red edge. Consider five vertices of $X$ that span both of these edges. The $K_5$ they form contains a blue and a red edge, so it is neither a red-green $K_5$ nor a blue-green $K_5$. By Theorem~\ref{thm:34n}, with $n=5$, it contains a blue-green $K_4$ and, applying the theorem again with blue and red interchanged, a red-green $K_4$.
\end{proof}

We also need the following result of Budden and DeJonge.

\begin{theorem}[{\cite[Theorem~3.5]{bd}}]\label{thm:355}
$g(3,5,5)=6$ and $g(3,5,6)=7$.
\end{theorem}

We note that for $6\le \min \{m,n\} \le9$, the following lemma is a consequence of \cite[Theorems~3.6--3.9]{bd}. However, we prove a more general result here.

\begin{figure}[htbp]
\centering
\begin{tikzpicture}[every node/.style={font=\small}]
 \begin{scope}[xshift=-1.45cm]
  \GRregular{a}{4}{.75}{135}
  \draw[rededge] (a0)--(a1) (a2)--(a3);
  \foreach \i/\j in {0/2,0/3,1/2,1/3}
    \draw[greenedge] (a\i)--(a\j);
  \foreach \i in {0,...,3} \node[vertex] at (a\i) {};
  \draw[black!65] (0,0) circle[radius=1.1];
  \node at (0,1.38) {$A$};
 \end{scope}

 \begin{scope}[xshift=1.45cm]
  \GRregular{b}{3}{.75}{90}
  \draw[blueedge] (b0)--(b1)--(b2);
  \draw[greenedge] (b0)--(b2);
  \foreach \i in {0,1,2} \node[vertex] at (b\i) {};
  \draw[black!65] (0,0) circle[radius=1.1];
  \node at (0,1.38) {$B$};
 \end{scope}

 \draw[greenedge] (-.35,0)--(.35,0);

\end{tikzpicture}
\caption{A partitioned $K_7$ with $m=n=6$ and no blue-red $K_3$. Blue edges are solid, red edges are dashed, and green edges are dotted. The vertices in $A$, together with the green edge in $B$, form a red-green $K_6$.}
\label{fig:3xy}
\end{figure}
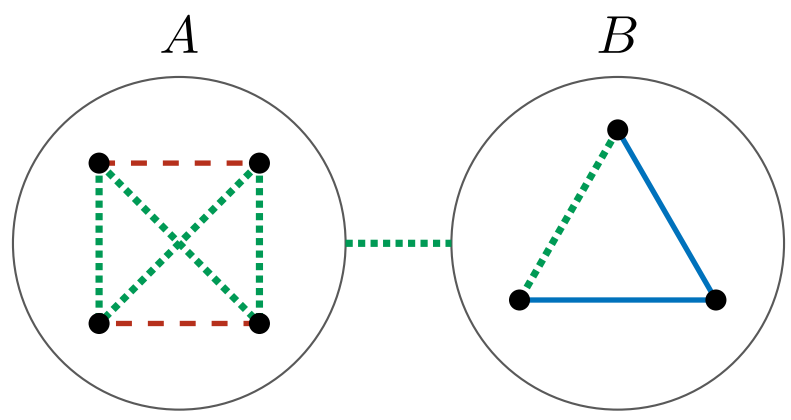

\begin{lemma}\label{thm:3xy}
For $m, n\ge6$, $g(3,m, n)\le m+n-5$. In particular,
$g(3,6,6)\le7$ and $g(3,6,7)\le8$.
\end{lemma}
\begin{proof}
Consider a Gallai three-coloring of $K_{m+n-5}$ containing no blue-red $K_3$. By Lemma~\ref{lem:tri2}, we can partition its vertices into the sets $A$ and $B$ described there. Since $|A|+|B|=m+n-5$, either $|A|\ge n-2$ or $|B|\ge m-2$.

If $|A|\ge n-2$, we show that there is a red-green $K_n$. Indeed, if $|A|\ge n$, this is immediate, since $A$ contains no blue edges. If $|A|=n-1$, then $|B|=m-4\ge2$, so $B\ne\emptyset$. Thus, $A$ together with a vertex from $B$ forms a red-green $K_n$. If $|A|=n-2$, then $|B|=m-3\ge3$, so $B$ contains a green edge to avoid a blue triangle. Combining $A$ with the endpoints of this green edge gives a red-green $K_n$ (see Figure~\ref{fig:3xy}).

A similar argument applies if $|B|\ge m-2$, giving a blue-green $K_m$.
\end{proof}

We again restate a theorem of Budden and DeJonge~\cite{bd}.

\begin{theorem}[{\cite[Theorem~4.1]{bd}}]\label{thm:44n}
For $n\ge4$,
\[
 g(4,4,n)=
 \begin{cases}
 n+1,&4\le n\le8,\\
 \min\{n,18\},&n\ge9.
 \end{cases}
\]
\end{theorem}

Here, we explicitly construct the colorings that show the lower bound for $4\le n\le8$, since the coloring $E_7$ will be used later for the lower bound of $g(6, 6, 7)$ in Theorem~\ref{thm:667}.

For $4\le n\le8$, define the following coloring $E_n$ of $K_n$. Consider two vertices $u$ and $v$, and color the edge $uv$ blue. The remaining $n-2$ vertices are partitioned into sets $P=\{p_1,p_2,\ldots\}$ and $Q=\{q_1,q_2,\ldots\}$ such that $|P|=|Q|$ if $n$ is even, and $|P|=|Q|+1$ otherwise. Use the following color scheme:
\begin{enumerate}
    \item Edges from $\{u,v\}$ to $P$ are red.
    \item Edges from $\{u,v\}$ to $Q$ are green.
    \item Edges inside $P$ are green.
    \item Edges inside $Q$ are red.
    \item Edges $p_iq_j$ between $P$ and $Q$ are red if $i+j$ is odd, and green otherwise.
\end{enumerate}
As an example, Figure~\ref{fig:En} shows $E_8$.

\begin{figure}[htbp]
\centering
\begin{tikzpicture}[
 every node/.style={font=\small},
 enregion/.style={
   draw=black!65,circle,minimum size=2.2cm,
   inner sep=0pt,outer sep=0pt
 }
]
 \node[enregion] (U) at (0,1.6) {};
 \node[enregion] (P) at (-1.85,-1.1) {};
 \node[enregion] (Q) at (1.85,-1.1) {};

 \node[above=2pt] at (U.north) {$\{u,v\}$};
 \node[above=2pt] at (P.north) {$P$};
 \node[above=2pt] at (Q.north) {$Q$};

 \draw[rededge] (U)--(P);
 \draw[greenedge] (U)--(Q);
 \draw[allowededge] (P)--(Q);

 \coordinate (u) at ([xshift=-.75cm]U.center);
 \coordinate (v) at ([xshift=.75cm]U.center);
 \draw[blueedge] (u)--(v);
 \node[vertex] at (u) {};
 \node[vertex] at (v) {};

 \foreach \prefix/\region/\edgestyle in
   {p/P/greenedge,q/Q/rededge} {
  \begin{scope}[shift={(\region.center)}]
   \foreach \i/\angle in {1/90,2/210,3/330}
    \coordinate (\prefix\i) at (\angle:.75);

   \draw[\edgestyle]
     (\prefix1)--(\prefix2)--(\prefix3)--cycle;

   \foreach \i in {1,2,3}
    \node[vertex] at (\prefix\i) {};
  \end{scope}
 }

 \begin{scope}[xshift=5.2cm,yshift=-1.25cm,text=black]
  \node at (1.05,2.95) {Coloring matrix for $p_i q_j$};

  \foreach \j in {1,2,3}
   \node at ({.7*(\j-.5)},2.35) {$q_\j$};

  \foreach \i in {1,2,3} {
   \node[anchor=east]
     at (-.15,{.7*(3.5-\i)}) {$p_\i$};

   \foreach \j in {1,2,3} {
    \pgfmathtruncatemacro{\parity}{mod(\i+\j,2)}
    \ifodd\parity
     \node at ({.7*(\j-.5)},{.7*(3.5-\i)}) {R};
    \else
     \node at ({.7*(\j-.5)},{.7*(3.5-\i)}) {G};
    \fi
   }
  }

  \draw[step=.7,black] (0,0) grid (2.1,2.1);
 \end{scope}

 \node at (2.5,-2.75)
   {$\wbr(E_8)=\wbg(E_8)=3,\qquad \wrg(E_8)=7$};

\end{tikzpicture}
\caption{The coloring $E_8$ in Theorem~\ref{thm:44n}. The hollow line between $P$ and $Q$ represents the nine edges between $P$ and $Q$, each red or green as specified by the matrix.}
\label{fig:En}
\end{figure}

Later, we will make use of the following lemma.

\begin{lemma}\label{lem:smallblocks}
Let $X$ be a Gallai three-colored complete graph with no red-green triangle. The following hold.
\begin{enumerate}[label=(\arabic*), ref=\thetheorem(\arabic*)]
\item\label{lem:sb1} If $|X|\ge3$, then $\wbr,\wbg\ge2$. If moreover $\wbr=2$, then $X$ is blue-green and $\wbg=|X|$; symmetrically, if $\wbg=2$, then $X$ is blue-red and $\wbr=|X|$.
\item\label{lem:sb2} If $|X|=4$ and $\wbr,\wbg\le3$, then $\wbr=\wbg=3$.
\item\label{lem:sb3} If $\wbr,\wbg\le4$, then $|X|\le5$. If moreover $|X|=5$, then $\wbr=\wbg=4$.
\item\label{lem:sb4} If $\wbr,\wbg\le5$, then $|X|\le6$. If moreover $|X|=6$, then $\wbr,\wbg\ge4$, and at least one of them equals $5$.
\end{enumerate}
\end{lemma}
\begin{proof}
For (1), $\wbr=1$ would make $X$ all green, and any three of its vertices would form a red-green triangle. If $\wbr=2$, then $X$ has neither a blue-red nor a red-green triangle, so Lemma~\ref{lem:tri3}, with blue and red interchanged, shows that $X$ has no red edge. The case $\wbg=2$ is symmetric.

For (2), by (1), both $\wbr$ and $\wbg$ are at least two, and if either equals two, then the other equals four. Hence $\wbr=\wbg=3$.

For (3), Theorem~\ref{thm:355} (with the colors permuted) gives $|X|\le5$. A graph on five vertices with $\wbr,\wbg\le4$ and $\wrg\le2$ uses all three colors, so permuting the colors in Corollary~\ref{cor:both4} gives $\wbr=\wbg=4$.

For (4), $g(3,6,6)\le7$ (Lemma~\ref{thm:3xy}, with the colors permuted) gives $|X|\le6$. If $|X|=6$, then $g(3,4,6)=6$ (Theorem~\ref{thm:34n}) forces each of $\wbr$ and $\wbg$ to be at least four, and Theorem~\ref{thm:355} forces one of them to be five (in both cases with the colors permuted).
\end{proof}

We also use two further values of Budden and DeJonge.

\begin{theorem}[{\cite[Theorem~4.2]{bd}}]\label{thm:455}
$g(4,5,5)=7$.
\end{theorem}

\begin{theorem}[{\cite[Theorem~4.2]{bd}}]\label{thm:456}
$g(4,5,6)=8$.
\end{theorem}

We again show the coloring that proves the lower bound in Theorem~\ref{thm:456}, since it will be useful later. It appears in Figure~\ref{fig:456lower}; we refer to its red-green $K_5$, whose red edges form a five-cycle, as $C$.

\begin{figure}[htbp]
\centering
\begin{tikzpicture}[
 every node/.style={font=\small},
 counterblock/.style={
  draw=black!65,circle,minimum size=2.2cm,
  inner sep=0pt,outer sep=0pt
 }
]
 \node[counterblock] (C) at (-1.75,0) {};
 \node[counterblock] (T) at (1.75,0) {};
 \draw[blueedge] (C)--(T);

 \begin{scope}[shift={(C.center)}]
  \foreach \i in {0,...,4}
   \coordinate (c\i) at ({90-72*\i}:.75);

  \foreach \i/\j in {0/2,0/3,1/3,1/4,2/4}
   \draw[greenedge] (c\i)--(c\j);
  \draw[rededge] (c0)--(c1)--(c2)--(c3)--(c4)--cycle;

  \foreach \i in {0,...,4}
   \node[vertex] at (c\i) {};
 \end{scope}

 \begin{scope}[shift={(T.center)}]
  \coordinate (t0) at (-.75,0);
  \coordinate (t1) at (.75,0);
  \draw[greenedge] (t0)--(t1);
  \node[vertex] at (t0) {};
  \node[vertex] at (t1) {};
 \end{scope}

 \node[above=4pt] at (C.north) {$C$};
 \node at (0,-1.65)
  {$\wbr=3,\qquad \wbg=4,\qquad \wrg=5$};
\end{tikzpicture}
\caption{A Gallai three-coloring of $K_7$ showing $g(4,5,6)\ge8$.}
\label{fig:456lower}
\end{figure}

\begin{theorem}\label{thm:466}
$g(4,6,6)=9$.
\end{theorem}
\begin{proof}
Suppose a Gallai three-coloring of $K_9$ avoids the three targets. Choose a minimal decomposition with base graph $H$ and $r$ blocks. Choosing one vertex from each of six blocks would give a two-colored $K_6$, so $r\le5$.

Suppose first that $r\ge4$. The base graph cannot be blue-red, since choosing one vertex from each of four blocks would give a blue-red $K_4$. Thus, by interchanging blue and red, we may assume without loss of generality that the base graph is blue-green.

Every block then has both a blue and a green neighbor. Consequently, every all-red block has order at most two; otherwise, three vertices of that block and one vertex of a blue neighbor would form a blue-red $K_4$.

Furthermore, there cannot be three all-red blocks of order two. If two of them are joined in blue, their four vertices form a blue-red $K_4$. Otherwise, the three blocks are pairwise joined in green, and their six vertices form a red-green $K_6$.

\case{Case 1: $r=5$}
Every edge within a block must be red. Indeed, a blue or green edge in one block, together with one vertex from each of the other four blocks, would give a blue-green $K_6$.

Each block therefore has order at most two, and at most two blocks can have order two. The total order is consequently at most $2+2+1+1+1=7$, which is a contradiction.

\case{Case 2: $r=4$}
At most one block contains a blue or green edge. Otherwise, take the endpoints of such an edge from each of two blocks and one vertex from each of the remaining two blocks. These six vertices form a blue-green $K_6$.

There must be a block that is not all red, since four all-red blocks would each have order at most two and could not contain all nine vertices. Call this block $X$. The other three blocks are all red. Each has order at most two, and they cannot all have order two, so their total order is at most five. Thus $|X|\ge4$.

There is no blue-green triangle in $X$, since it, along with one vertex from each of the other three blocks, would form a blue-green $K_6$. There is also no blue-red triangle in $X$, since adding a vertex from a blue neighbor would give a blue-red $K_4$. Lemma~\ref{lem:tri3} therefore implies that $X$ is red-green. If $|X|\ge5$, five vertices of $X$ and one vertex from a green neighbor would form a red-green $K_6$, so $|X|=4$. Hence, the other three blocks have orders $2$, $2$, and $1$.

Moreover, $X$ contains a red edge, since otherwise it would be all green and would contain a blue-green triangle. Choose an all-red block $D$ of order two. If $X$ and $D$ are joined in blue, a red edge in $X$ and the red edge in $D$ give a blue-red $K_4$. If they are joined in green, all four vertices of $X$ and both vertices of $D$ give a red-green $K_6$. Both possibilities are forbidden, as illustrated in Figure~\ref{fig:466-exception-step}.

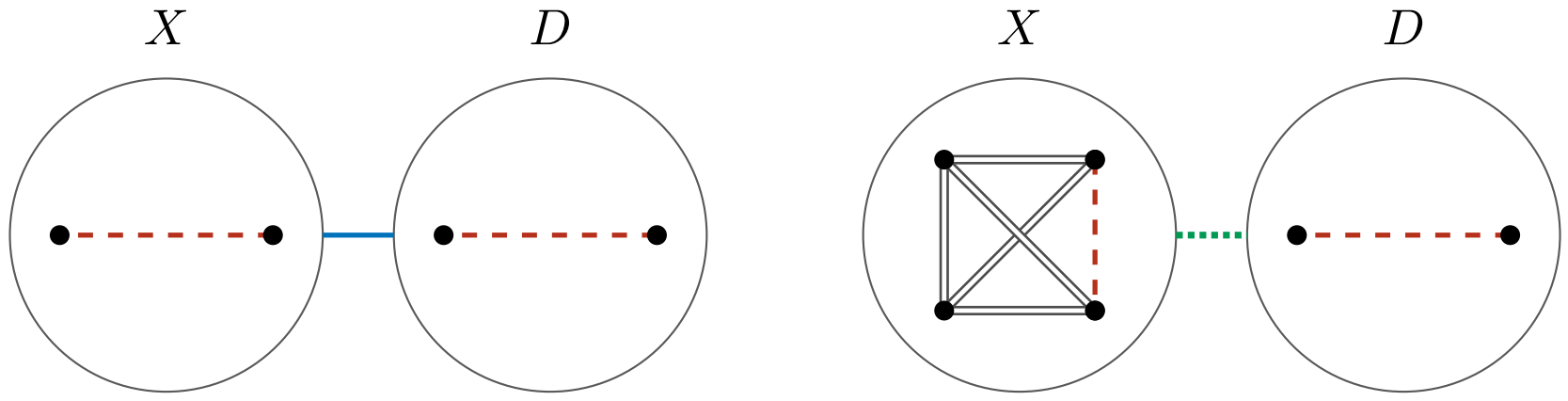
\begin{figure}[htbp]
\centering
\begin{tikzpicture}[every node/.style={font=\small}]

  \foreach \v/\x/\n/\ang/\sty/\lab/\txt/\nodraw in {
    Xb/0/2/0/rededge/X/{}/0,
    Db/2.7/2/0/rededge/D/{}/0,
    Xg/6/4/45/allowededge/X/{}/12,
    Dg/8.7/2/0/rededge/D/{}/0
  } {
    \node[draw=black!65,circle,minimum size=2.2cm,
          inner sep=0pt,outer sep=0pt] (\v) at (\x,0) {};
    \foreach \i in {1,...,\n} {
      \coordinate (\v-\i) at
        ($(\v.center)+({\ang-360*(\i-1)/\n}:.75)$);
    }
    \foreach \i in {1,...,\n} {
      \foreach \j in {1,...,\n} {
        \ifnum\i<\j\relax
          \ifnum\numexpr10*\i+\j\relax=\nodraw\relax\else
            \draw[\sty] (\v-\i)--(\v-\j);
          \fi
        \fi
      }
    }
    \foreach \i in {1,...,\n} {
      \node[vertex] at (\v-\i) {};
    }
    \node[above=3pt] at (\v.north) {$\lab$};
    \node[below=3pt,font=\scriptsize,align=center,
          text width=2.2cm] at (\v.south) {\txt};
  }
  \draw[rededge] (Xg-1)--(Xg-2);
  \foreach \i in {1,2} {
    \node[vertex] at (Xg-\i) {};
  }
  \draw[blueedge] (Xb)--(Db);
  \draw[greenedge] (Xg)--(Dg);
\end{tikzpicture}
\caption{$X$ is red-green and contains a red edge, and $D$ is an all-red block of order two. In the first panel, $X$ and $D$ are joined in blue; in the second, they are joined in green, and the remaining edges of $X$ are red or green.}
\label{fig:466-exception-step}
\end{figure}

\case{Case 3: $r=2$}
Let $A$ and $B$ be the two blocks.

\case{Case 3a: $A$ and $B$ are joined in green}
Neither block contains a blue-green or red-green $K_5$, since such a clique extends to a forbidden $K_6$ by adding a vertex of the other block. Both blocks also avoid a blue-red $K_4$. Theorem~\ref{thm:455} therefore bounds the order of each block by six.

Without loss of generality, suppose that $|A|\ge|B|$. Since $|A|+|B|=9$, the possible orders are $6+3$ and $5+4$.

Suppose first that $|A|=6$ and $|B|=3$. If $B$ is all red, then $A$ has no red-green triangle, since such a triangle together with $B$ would form a red-green $K_6$. Theorem~\ref{thm:355} (with the colors permuted) then gives a blue-green $K_5$ in $A$, which together with a vertex of $B$ forms a blue-green $K_6$. Otherwise, $B$ contains a blue or green edge. There can be no blue-green $K_4$ in $A$, since its vertices and the endpoints of this edge would form a blue-green $K_6$. Thus $A$ has no blue-red $K_4$, no blue-green $K_4$, and no red-green $K_5$. This contradicts $g(4,4,5)=6$ from Theorem~\ref{thm:44n}.

Now suppose that $|A|=5$ and $|B|=4$. The block $B$ contains a green edge, since otherwise it is a blue-red $K_4$. Since $\gr(K_4)=5$, the block $A$ contains a two-colored $K_4$. A blue-red $K_4$ is already forbidden. A blue-green or red-green $K_4$, together with the endpoints of the green edge in $B$, gives the corresponding forbidden $K_6$.

\case{Case 3b: $A$ and $B$ are joined in blue or red}
By interchanging blue and red, we may assume without loss of generality that they are joined in blue.

At least one block must be all green. Otherwise, a blue or red edge in each block would give a blue-red $K_4$. Relabeling if necessary, let $B$ be an all-green block.

If $|B|\ge5$, five vertices of $B$ and one vertex of $A$ give a blue-green $K_6$. Hence $|B|\le4$, so $|A|\ge5$.

There is no blue-red triangle in $A$, since adding a vertex of $B$ would give a blue-red $K_4$. Similarly, $A$ has no blue-green $K_5$, since adding a vertex of $B$ would give a blue-green $K_6$. It also has no red-green $K_6$ by assumption. The equality $g(3,5,6)=7$ from Theorem~\ref{thm:355} therefore gives $|A|\le6$. Since $|A|+|B|=9$, we have $|B|\ge3$.

If $|B|=3$, then $|A|=6$. A blue-green triangle in $A$, together with $B$, would form a blue-green $K_6$. Thus $A$ has neither a blue-red nor a blue-green triangle, so by Lemma~\ref{lem:tri3}, its six vertices form a red-green $K_6$, a contradiction.

If $|B|=4$, then $|A|=5$. Any blue or green edge in $A$, together with $B$, would form a blue-green $K_6$. Hence every edge of $A$ is red, giving a blue-red $K_4$.

This completes the upper bound. For equality, consider Figure~\ref{fig:466lower}. It consists of two blocks joined in blue: the first is the coloring $C$ from Theorem~\ref{thm:456}, and the second, $T$, is a green $K_3$.
\end{proof}

\begin{figure}[htbp]
\centering
\begin{tikzpicture}[
 every node/.style={font=\small},
 counterblock/.style={
  draw=black!65,circle,minimum size=2.2cm,
  inner sep=0pt,outer sep=0pt
 }
]
 \node[counterblock] (C) at (-1.75,0) {};
 \node[counterblock] (T) at (1.75,0) {};

 \draw[blueedge] (C)--(T);

 \begin{scope}[shift={(C.center)}]
  \foreach \i in {0,...,4}
   \coordinate (c\i) at ({90-72*\i}:.75);

  \foreach \i/\j in {0/2,0/3,1/3,1/4,2/4}
   \draw[greenedge] (c\i)--(c\j);

  \draw[rededge]
   (c0)--(c1)--(c2)--(c3)--(c4)--cycle;

  \foreach \i in {0,...,4}
   \node[vertex] at (c\i) {};
 \end{scope}

 \begin{scope}[shift={(T.center)}]
  \foreach \i/\angle in {0/90,1/210,2/330}
   \coordinate (t\i) at (\angle:.75);

  \draw[greenedge] (t0)--(t1)--(t2)--cycle;

  \foreach \i in {0,1,2}
   \node[vertex] at (t\i) {};
 \end{scope}

 \node[above=4pt] at (C.north) {$C$};
 \node[above=4pt] at (T.north) {$T$};

 \node at (0,-1.65)
  {$\wbr=3,\qquad \wbg=\wrg=5$};
\end{tikzpicture}
\caption{A Gallai three-coloring of $K_8$ showing $g(4,6,6)\ge9$.}
\label{fig:466lower}
\end{figure}

Next, we determine $g(5,5,7)$, which we use in Sections~\ref{sec:667} and~\ref{sec:k7}.

\begin{theorem}\label{thm:557}
$g(5,5,7)=11$.
\end{theorem}
\begin{proof}
Suppose a Gallai three-coloring of $K_{11}$ has no blue-red or blue-green $K_5$ and no red-green $K_7$. As before, choose a minimal decomposition with base graph $H$ and $r$ blocks. Since one vertex from each block forms a two-colored clique, $r\le6$.

\case{Case 1: the base graph is blue-red or blue-green, and $r\ge4$}
By symmetry, we may assume that the base graph is blue-red. Then $r=4$, and all blocks are all green, since we otherwise obtain a blue-red $K_5$. Every block has a blue neighbor to prevent a monochromatic cut, so its order is at most three (otherwise, that block and a vertex from the blue neighbor would give a blue-green $K_5$). To total eleven, the orders must be $3,3,3,2$. But any two blocks joined in blue then give a blue-green $K_5$.

\case{Case 2: the base graph is red-green, and $r=6$}
Every block must be all blue. If not, a red or green edge in one block and a vertex from every other block form a red-green $K_7$. Additionally, any two blocks have combined order at most four to prevent a forbidden $K_5$. Eleven vertices force the orders $2,2,2,2,2,1$. By $R(3,3)=6$, the base graph $H$ has a monochromatic triangle, and the corresponding blocks have total order at least five, giving a forbidden $K_5$.

\case{Case 3: the base graph is red-green, and $r=5$}
By similar arguments, there is at most one block $X$ that is not blue, and it has $\wrg(X)\le2$ and $\wbr(X),\wbg(X)\le3$. So $|X|\le4$, since $\gr(K_4) = 5$.

Any two blue blocks have combined order at most four to prevent a forbidden $K_5$. So, if all blocks are blue, then the total order is at most ten, a contradiction. Otherwise, the four blue blocks have total order at most eight, so $|X|\ge3$.

If $|X|=4$, Lemma~\ref{lem:sb2} gives $\wbr(X)=\wbg(X)=3$, which forces the other four blocks to be singletons, contradicting the total number of vertices. If $|X|=3$, the other four blocks have total order eight and pairwise combined order at most four, so each has order two. Since every other block is blue and $X$ has both a red and a green neighbor, we require $\wbr(X),\wbg(X)\le2$ to prevent the forbidden cliques. This is a contradiction, since $X$ is a Gallai triangle with $\wrg(X)\le2$.

\case{Case 4: the base graph is red-green, and $r=4$}
All blocks must satisfy $\wbr,\wbg\le3$. As above, any two blue blocks have combined order at most four. At most two blocks contain a red or green edge, since three such edges and a vertex from the fourth block form a red-green $K_7$.

If all the blocks are blue, their total order is at most eight, a contradiction.

If three are blue, call the remaining block $X$. It must have $\wrg(X)\le3$, since a red-green $K_4$ in $X$, together with one vertex from each other block, would form a red-green $K_7$. Since $\gr(K_4)=5$, we get $|X| \leq 4$. The other three blocks have pairwise combined order at most four, so they have total order at most six. Hence the total order is at most ten, a contradiction.

If there are two blocks $X$ and $Y$ that are not blue, both must have $\wrg\le2$, since a red-green triangle in one of them, a red or green edge in the other, and one vertex from each blue block would form a red-green $K_7$. Since $\gr(K_4)=5$, both have order at most four. The two blue blocks together have order at most four. Since there are eleven vertices, the orders of $X$ and $Y$ must be $4+4$ or $4+3$. Assume $|X| \geq |Y|$. Then, by Lemma~\ref{lem:sb2}, $X$ has $\wbr=\wbg=3$, and by Lemma~\ref{lem:sb1}, $Y$ has $\wbr,\wbg\ge2$. Whichever color joins $X$ and $Y$, this results in a forbidden $K_5$, as illustrated in Figure~\ref{fig:557-exceptions-step}.

\begin{figure}[htbp]
\centering
\begin{tikzpicture}[every node/.style={font=\small}]

  \foreach \v/\x/\n/\ang/\sty/\lab/\txt in {
    Xr/0/3/90/allowededge/X/{},
    Yr/2.7/2/0/allowededge/Y/{},
    Xg/6/3/90/allowededge/X/{},
    Yg/8.7/2/0/allowededge/Y/{}
  } {
    \node[draw=black!65,circle,minimum size=2.2cm,
          inner sep=0pt,outer sep=0pt] (\v) at (\x,0) {};
    \foreach \i in {1,...,\n} {
      \coordinate (\v-\i) at
        ($(\v.center)+({\ang-360*(\i-1)/\n}:.75)$);
    }
    \foreach \i in {1,...,\n} {
      \foreach \j in {1,...,\n} {
        \ifnum\i<\j\relax
          \draw[\sty] (\v-\i)--(\v-\j);
        \fi
      }
    }
    \foreach \i in {1,...,\n} {
      \node[vertex] at (\v-\i) {};
    }
    \node[above=3pt] at (\v.north) {$\lab$};
    \node[below=3pt,font=\scriptsize,align=center,
          text width=2.2cm] at (\v.south) {\txt};
  }
  \draw[rededge] (Xr)--(Yr);
  \draw[greenedge] (Xg)--(Yg);
\end{tikzpicture}
\caption{$|X|=4$ and $|Y|\ge3$. In the first panel, $X$ contains a blue-red $K_3$ and $Y$ contains a blue or red edge. In the second panel, $X$ contains a blue-green $K_3$ and $Y$ contains a blue or green edge.}
\label{fig:557-exceptions-step}
\end{figure}
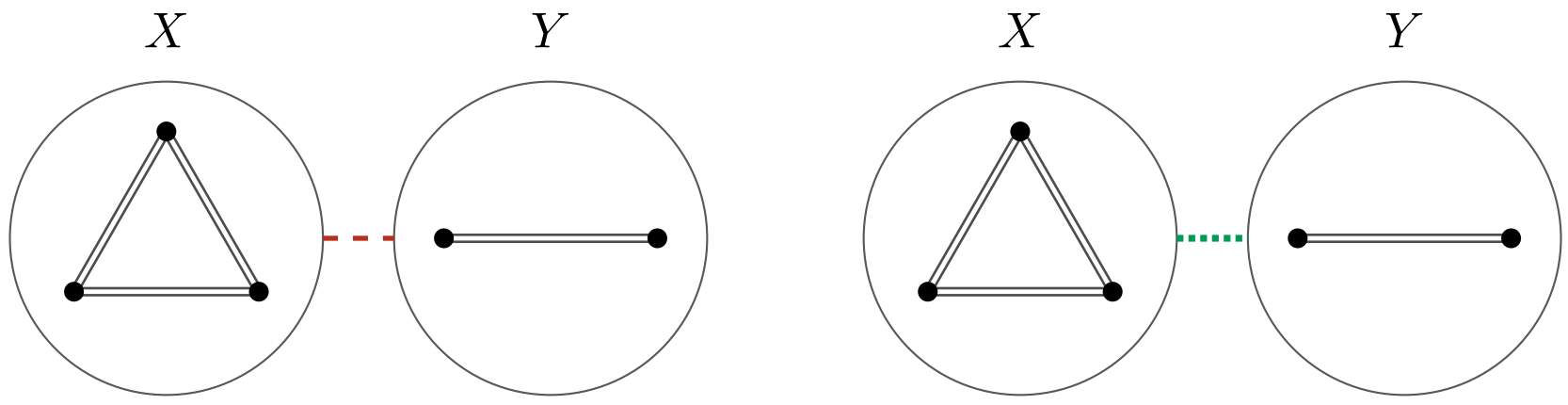

\case{Case 5: $r=2$}
Let $A$ and $B$ be the two blocks, with $|A|\ge|B|$. We consider the three possible colors joining them.

If they are joined in blue, both blocks have $\wbr,\wbg\le3$ and $\wrg\le6$. Since $g(4,4,7)=8$ (Theorem~\ref{thm:44n}), we have $|A|,|B|\le7$. Hence $|A|\ge6$ and $|B|\ge4$. If $\wbr(B)=1$, then $B$ is all green and, together with one vertex of $A$, forms a blue-green $K_5$; similarly, $\wbg(B)=1$ would make $B$ all red and give a blue-red $K_5$. So $B$ contains a blue or red edge and a blue or green edge, and hence $\wbr(A),\wbg(A)\le2$. By Lemma~\ref{lem:tri3}, $A$ is red-green, and $R(3,3)=6$ gives a red or green triangle in $A$, contradicting $\wbr(A),\wbg(A)\le2$.

If they are joined in red, both blocks have $\wbr\le3$, $\wbg\le4$, and $\wrg\le5$. So, by Theorem~\ref{thm:456}, $7 \geq |A| \geq |B|$. Hence $|A|\ge6$ and $|B|\ge4$. If $\wbr(B)=1$, then $B$ is all green, which forces $\wrg(A)\le2$ to prevent a red-green $K_7$. However, together with $\wbr(A)\le3$ and $\wbg(A)\le4$, this contradicts $g(3,4,5)=5$ (Theorem~\ref{thm:34n}, with the colors permuted). Hence $\wbr(B)\ge2$, and so $\wbr(A)\le2$. Theorem~\ref{thm:355} applied to $A$ then gives $\wrg(A)\ge5$, since $\wbg(A)\le4$. To prevent a red-green $K_7$, this forces $B$ to be blue. This, however, yields a blue-red $K_5$.

The case where they are joined in green follows similarly by interchanging red and green. For equality, see Figure~\ref{fig:557lower}: both blocks contain the coloring $C$ from Theorem~\ref{thm:456}. This is, in fact, the $n=6$ coloring from \cite[Theorem~2.4]{bw} with permuted colors.
\end{proof}

\begin{figure}[htbp]
\centering
\begin{tikzpicture}[every node/.style={font=\small}]
  \node[draw=black!65,circle,minimum size=2.2cm,
        inner sep=0pt,outer sep=0pt] (A) at (0,0) {};
  \node[draw=black!65,circle,minimum size=2.2cm,
        inner sep=0pt,outer sep=0pt] (B) at (3.5,0) {};

  \draw[blueedge] (A)--(B);

  \foreach \X in {A,B} {
    \foreach \i in {0,...,4} {
      \coordinate (\X-\i) at
        ($(\X.center)+({90-72*\i}:.75)$);
    }

    \foreach \i/\j in {0/1,1/2,2/3,3/4,4/0} {
      \draw[rededge] (\X-\i)--(\X-\j);
    }

    \foreach \i/\j in {0/2,0/3,1/3,1/4,2/4} {
      \draw[greenedge] (\X-\i)--(\X-\j);
    }

    \foreach \i in {0,...,4} {
      \node[vertex] at (\X-\i) {};
    }
    \node[above=4pt] at (\X.north) {$C$};
  }
\end{tikzpicture}
\caption{A Gallai three-coloring of $K_{10}$ showing
$g(5,5,7)\ge11$.}
\label{fig:557lower}
\end{figure}

\section{The off-diagonal value \texorpdfstring{$g(6,6,7)$}{g(6,6,7)}} \label{sec:667}

In this section, we determine $g(6,6,7)$, which handles the two-block case in the proof of Theorem~\ref{thm:k7}.

\begin{theorem}\label{thm:667}
$g(6,6,7)=13$.
\end{theorem}

\begin{proof}
Consider a Gallai three-coloring of $K_{13}$ that contains no blue-red $K_6$, no blue-green $K_6$, and no red-green $K_7$. Choose a minimal decomposition, with base graph $H$ of order $r$. Choosing one vertex from each block shows that $r\le6$. Thus $r=2$, $4$, $5$, or $6$.

As noted in Section~\ref{sec:prelim}, when $r\ge4$, every block has a neighbor in each base color. Additionally, when $r=4$, every triangle of the base graph uses both colors. If not, a monochromatic triangle would force all three edges to the remaining vertex to have the other color, giving a monochromatic cut and contradicting minimality.

We may interchange red and green by symmetry.

\case{Case 1: $r=6$}
The base colors must be red and green; otherwise, one vertex from each block gives a forbidden $K_6$. All six blocks must be blue; otherwise, we obtain a red-green $K_7$.

Any two blocks have combined order at most five to prevent a forbidden $K_6$. One can verify that this implies the only possibility is $(3,2,2,2,2,2)$. By $R(3,3)=6$, the base graph contains a monochromatic triangle. The corresponding blocks contain at least six vertices and give a forbidden $K_6$.

\case{Case 2: $r=5$ and the base colors include blue}
Assume the base colors to be blue and red. A blue or red edge inside a block, together with one vertex from each of the other four blocks, would give a blue-red $K_6$. Thus every block is all green.

Any two blocks have combined order at most six, since together they form a blue-green or red-green clique. Hence all blocks have order at most three, since a block of order at least four would leave at most two vertices in each other block, for a total of at most twelve. So, there are at least three blocks of order three. Two of these joined in blue would give a blue-green $K_6$. Hence these three blocks are pairwise joined in red, giving a red-green $K_9$, again a contradiction.

\case{Case 3: $r=5$ and the base colors are red and green}
No block contains a red-green triangle, since adding one vertex from each of the other four blocks would give a red-green $K_7$. Moreover, at most one block contains a red or green edge, as two such edges and one vertex from each remaining block would also give a red-green $K_7$.

Any two all-blue blocks have combined order at most five, so at most one of them has order greater than two, and if it has order $k\ge3$, every other all-blue block has order at most $5-k\le2$. Hence five all-blue blocks have total order at most $3+2+2+2+2=11$, and four have total order at most $3+2+2+2=9$.

There must therefore be a block $X$ that is not all blue, and $|X|\ge4$. A blue-red or blue-green $K_5$ in $X$ extends by one vertex of a red or green neighbor, respectively. Thus neither occurs in $X$, and Lemma~\ref{lem:sb3} implies $|X|\le5$.

If $|X|=5$, Lemma~\ref{lem:sb3} gives both a blue-red $K_4$ and a blue-green $K_4$ in $X$. Any other block (which is all blue) with at least two vertices would extend one of these to a forbidden $K_6$. Hence, the other four blocks are singletons, which is a contradiction.

If $|X|=4$, the other orders must be $3,2,2,2$. Suppose that $X$ is joined in red to the block of order three (the other case is again symmetric). Then $X$ has no blue-red triangle. Since it also has no red-green triangle, Lemma~\ref{lem:sb1} implies that $X$ is entirely blue-green. Since the base graph has no monochromatic cut, $X$ has a green neighbor $Z$, which must be an all-blue block of order at least two. Then $X$ and two vertices of $Z$ give a blue-green $K_6$ (see Figure~\ref{fig:exceptionsteps}).

\begin{figure}[htbp]
\centering
\begin{tikzpicture}[every node/.style={font=\small}]
  \foreach \v/\x/\n/\ang/\sty/\lab/\txt in {
    Xg/6/4/45/allowededge/X/{},
    Zg/8.7/2/0/blueedge/Z/{}
  } {
    \node[draw=black!65,circle,minimum size=2.2cm,
          inner sep=0pt,outer sep=0pt] (\v) at (\x,0) {};
    \foreach \i in {1,...,\n} {
      \coordinate (\v-\i) at
        ($(\v.center)+({\ang-360*(\i-1)/\n}:.75)$);
    }
    \foreach \i in {1,...,\n} {
      \foreach \j in {1,...,\n} {
        \ifnum\i<\j\relax
          \draw[\sty] (\v-\i)--(\v-\j);
        \fi
      }
    }
    \foreach \i in {1,...,\n} {
      \node[vertex] at (\v-\i) {};
    }
    \node[above=3pt] at (\v.north) {$\lab$};
    \node[below=3pt,font=\scriptsize,align=center,
          text width=2.2cm] at (\v.south) {\txt};
  }
  \draw[greenedge] (Xg)--(Zg);
\end{tikzpicture}
\caption{$X$ is a blue-green $K_4$, and $Z$ contains a blue edge.}
\label{fig:exceptionsteps}
\end{figure}
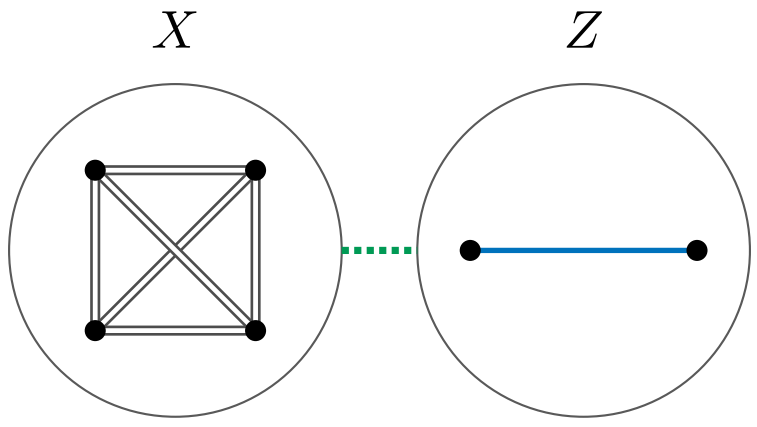

\case{Case 4: $r=4$ and the base colors include blue}
Again, take the base colors to be blue and red. At most one block contains a blue or red edge, since two such edges and one vertex from each of the other two blocks would give a blue-red $K_6$. If all four blocks are green, each pair has combined order at most six, and hence the total order is at most twelve, which is a contradiction. Hence there is one block $X$ that is not all green.

The block $X$ contains no blue-red triangle, since such a triangle extends by one vertex from each other block to a blue-red $K_6$. It also contains no blue-green $K_5$ and no red-green $K_6$, because it has both a blue neighbor and a red neighbor. The equality $g(3,5,6)=7$ from Theorem~\ref{thm:355} thus implies $|X|\le6$.

The three remaining blocks are green. Their base triangle uses both blue and red (because $r=4$). Two of these blocks joined in blue have combined order at most five, and both other pairs have combined order at most six. Adding these three bounds counts each vertex twice, so the three blocks have total order at most eight. Thus, $|X|\ge5$.

If $X$ uses all three colors, then by Corollary~\ref{cor:both4}, it has both a blue-green $K_4$ and a red-green $K_4$. Each remaining green block then has order at most two, and hence the total order is at most $6+2+2+2=12$.

So, $X$ uses at most two colors. It is not blue-red, since it has no blue-red triangle, and it is not blue-green, since it has no blue-green $K_5$. Hence, $X$ is red-green. Since $X$ has no red-green $K_6$, we have $|X|=5$. Every red neighbor of $X$ is consequently a singleton. Also, $X$ contains a green edge, since it has no blue-red triangle. Hence, every blue neighbor of $X$ has order at most three, since it is an all-green block. Since there is at least one red neighbor, the total order is at most $5+1+3+3=12$, a contradiction.

\case{Case 5: $r=4$ and the base colors are red and green}
Each block has a red and a green neighbor, so no block contains a blue-red or blue-green $K_5$. At most two blocks contain a red or green edge: three such edges and one vertex from the fourth block would give a red-green $K_7$.

If every block is blue, each pair has combined order at most five, and the total is at most ten.

Suppose exactly one block $X$ is not all blue. It contains no red-green $K_4$, since one vertex from each other block would yield a red-green $K_7$. Theorem~\ref{thm:455}, with the red-green target of order four, therefore gives $|X|\le6$.

The three blue blocks have total order at most seven (by pairwise bounds of five), so $|X|=6$. We recall $g(4,4,5)=6$ from Theorem~\ref{thm:44n}, and permute the colors. Hence, $X$ has both a blue-red $K_4$ and a blue-green $K_4$. A similar argument to the one used in Case~3 ensures that every remaining block is a singleton, giving only nine vertices.

It remains to consider two blocks $X$ and $Y$ that are not all blue. Neither contains a red-green triangle, to prevent a red-green $K_7$. Lemma~\ref{lem:sb3} gives $|X|,|Y|\le5$. The two blue blocks have combined order at most five, as above. Assuming $|X|\ge|Y|$, thirteen vertices therefore require either $|X|=5$ and $|Y|\ge3$, or $|X|=|Y|=4$.

In the first case, Lemma~\ref{lem:sb3} gives both a blue-red $K_4$ and a blue-green $K_4$ in $X$, while Lemma~\ref{lem:sb1} gives an edge in each of those color pairs in $Y$. Choosing the pair containing the color between $X$ and $Y$ gives a forbidden $K_6$.

Now let $|X|=|Y|=4$. By symmetry, we may assume that $X$ and $Y$ are joined in red. They cannot both contain a blue-red triangle; say $X$ contains none. Lemma~\ref{lem:sb1} then makes $X$ a blue-green $K_4$. The green neighbor of $X$ must be an all-blue block (because $X$ and $Y$ are joined in red), and it must be a singleton, since otherwise $X$ and two of its vertices would form a blue-green $K_6$. Since $X$ is blue-green and has no red-green triangle, it contains a blue edge. This edge and four vertices of an all-blue block would give a forbidden $K_6$, so the remaining blue block has order at most three. Hence, the total order is at most $4+4+1+3=12$.

\case{Case 6: $r=2$ and the two blocks are joined in blue}
Let $A$ and $B$ be the two blocks, with $|A|\ge|B|$, and write $a=|A|$ and $b=|B|$. Neither block contains a blue-red or blue-green $K_5$, since either extends by one vertex of the other block. Neither contains a red-green $K_7$. Theorem~\ref{thm:557} therefore gives $a,b\le10$, so $b\ge3$.

Suppose first that $a\ge8$. By $g(4,4,7)=8$ from Theorem~\ref{thm:44n}, $A$ contains a blue-red or a blue-green $K_4$. By symmetry, we may assume that it is blue-red. Any blue or red edge in $B$ would extend it to a blue-red $K_6$, so $B$ is all green. Since $b\ge3$, a blue-green triangle in $A$ would give a blue-green $K_6$, so $A$ has no blue-green triangle. Then, $g(3,6,7)\le8$ from Lemma~\ref{thm:3xy}, with the colors permuted, implies that $A$ has a blue-red $K_6$ or a red-green $K_7$, a contradiction.

The remaining orders are $a=7$ and $b=6$. A blue-red $K_4$ in $A$ would force $B$ to be an all-green $K_6$, and a blue-green $K_4$ would force an all-red $K_6$. Thus neither kind of $K_4$ occurs in $A$. Applying $g(3,4,7)=7$ from Theorem~\ref{thm:34n}, and then interchanging red and green, we get a blue-red triangle and a blue-green triangle in $A$. Consequently, $B$ contains neither of these triangles. Since $B$ is Gallai, all its triangles must therefore be red-green, so $B$ has no blue edge. But $R(3,3)=6$ gives a red or green triangle in $B$, a contradiction.

\case{Case 7: $r=2$ and the two blocks are joined in red}
The green case follows symmetrically. Again, let $A$ and $B$ be the two blocks, with $|A|\ge|B|$, and write $a=|A|$ and $b=|B|$. Neither block contains a blue-red $K_5$, a blue-green $K_6$, or a red-green $K_6$. Thus neither contains a two-colored $K_6$. Since $\gr(K_6)=12$, we get $a,b\le11$, and so $b\ge2$.

If $A$ contains a blue-red $K_4$, then $B$ is all green. Hence, in such a case, $b\le5$ and $a\ge8$. Now, if $b=2$, then $a=11$, and $A$ must avoid a blue-red $K_5$, a red-green $K_5$, and a blue-green $K_6$. This contradicts Theorem~\ref{thm:557} with blue and red interchanged. If $b\ge3$, then $A$ must avoid a red-green $K_4$, a blue-red $K_5$, and a blue-green $K_6$, which contradicts Theorem~\ref{thm:456} (with the colors permuted), since $a\ge8$.

So $A$ has no blue-red $K_4$. Theorem~\ref{thm:466} then implies $a\le8$. If $a=8$, then $b=5$, and Theorem~\ref{thm:456}, with the colors permuted, implies that $A$ contains a red-green $K_5$. Hence $B$ has to be all blue. Its five vertices and one vertex of $A$ then give a blue-red $K_6$.

The only remaining case is $a=7$ and $b=6$. The bound $g(3,6,6)\le7$ from Lemma~\ref{thm:3xy} forces a blue-red triangle in $A$, so $B$ contains no such triangle. Applying $g(3,4,6)=6$ from Theorem~\ref{thm:34n} gives a red-green $K_4$ in $B$. Hence $A$ has no red-green triangle. The same bound $g(3,6,6)\le7$, with the colors permuted, gives a contradiction, since $A$ has no blue-red or blue-green $K_6$.

This proves the upper bound. For equality, join the coloring $C$ (from Theorem~\ref{thm:456}) to the coloring $E_7$ (from Theorem~\ref{thm:44n}) in blue. This yields Figure~\ref{fig:667lower}.
\end{proof}

\begin{figure}[htbp]
\centering
\begin{tikzpicture}[every node/.style={font=\small}]
  \node[draw=black!65,circle,minimum size=2.2cm,
        inner sep=0pt,outer sep=0pt] (C) at (0,0) {};
  \node[draw=black!65,circle,minimum size=2.2cm,
        inner sep=0pt,outer sep=0pt] (E) at (3.5,0) {};
  \draw[blueedge] (C)--(E);

  \foreach \i in {0,...,4} {
    \coordinate (c\i) at
      ($(C.center)+({90-72*\i}:.75)$);
  }
  \foreach \i/\j in {0/1,1/2,2/3,3/4,4/0} {
    \draw[rededge] (c\i)--(c\j);
  }
  \foreach \i/\j in {0/2,0/3,1/3,1/4,2/4} {
    \draw[greenedge] (c\i)--(c\j);
  }
  \foreach \i in {0,...,4} {
    \node[vertex] at (c\i) {};
  }

  \foreach \i in {0,...,6} {
    \coordinate (e\i) at
      ($(E.center)+({90-360*\i/7}:.9)$);
  }
  \draw[blueedge] (e0)--(e1);
  \foreach \i/\j in
    {0/2,0/3,0/4,1/2,1/3,1/4,5/6,2/6,3/5,4/6} {
    \draw[rededge] (e\i)--(e\j);
  }
  \foreach \i/\j in
    {2/3,2/4,3/4,0/5,0/6,1/5,1/6,2/5,3/6,4/5} {
    \draw[greenedge] (e\i)--(e\j);
  }
  \foreach \i in {0,...,6} {
    \node[vertex] at (e\i) {};
  }

  \node[above=4pt] at (C.north) {$C$};
  \node[above=4pt] at (E.north) {$E_7$};
\end{tikzpicture}
\caption{A Gallai three-coloring of $K_{12}$ showing
$g(6,6,7)\ge13$.}
\label{fig:667lower}
\end{figure}

\section{The diagonal value on seven vertices}
\label{sec:k7}

We now combine the results above to prove our main result.

\begin{theorem}\label{thm:k7}
$\gr(K_7)=17$.
\end{theorem}

\begin{proof}
We first prove the upper bound. Suppose that a Gallai three-coloring of $K_{17}$ contains no two-colored $K_7$, and choose a minimal decomposition with base graph $H$ of order $r$. Choosing one vertex from each block gives $r\le6$. Thus $r=2$, $4$, $5$, or $6$.

For $r\ge4$, by symmetry, let the base colors be red and green. Every block has a neighbor in each color, so a blue-red or blue-green $K_6$ inside a block would extend to a two-colored $K_7$. Hence, we exclude both throughout these cases.

\case{Case 1: $r=6$}
A red or green edge in a block, together with one vertex from each of the other five blocks, gives a red-green $K_7$. Hence all blocks are blue, and any two have combined order at most six. If the largest order is at least four, the total number of vertices is at most fourteen, a contradiction.

Hence, every order is at most three, and seventeen vertices force the distribution $(3,3,3,3,3,2)$. Since $R(3,3)=6$, the base graph contains a monochromatic triangle. The corresponding blue blocks contain at least eight vertices, giving a two-colored $K_7$.

\case{Case 2: $r=5$}
As in Case~3 of the proof of Theorem~\ref{thm:667}, at most one block is not all blue, and any such block $X$ contains no red-green triangle. Any two blue blocks have combined order at most six. If all five blocks are blue, adding the pairwise bounds then gives total order at most fifteen.

Otherwise, Lemma~\ref{lem:sb4} gives $|X|\le6$. The four blue blocks have total order at most twelve, implying $5\le|X|\le6$.

If $|X|=5$, the other four orders total twelve and have pairwise sums at most six, so they are all three. A blue-red or blue-green $K_4$ in $X$ would extend to a forbidden $K_7$ using a red or green neighbor, respectively. Hence, $X$ contains no blue-red or blue-green $K_4$. Since $X$ also contains no red-green triangle, it contains no two-colored $K_4$. This contradicts $\gr(K_4)=5$.

If $|X|=6$, Lemma~\ref{lem:sb4} gives both a blue-red $K_4$ and a blue-green $K_4$ in $X$. Every all-blue block must then have order at most two. This is illustrated in Figure~\ref{fig:k7-exceptions-step}. The total order is at most $6+4\cdot2=14$, which is again impossible.

\case{Case 3: $r=4$}
As in Case~5 of the proof of Theorem~\ref{thm:667}, at most two blocks are not all blue. If every block is blue, the total order is at most twelve.

Suppose exactly one block $X$ is not all blue. It has no red-green $K_4$, by the same argument as in Case~5 of the proof of Theorem~\ref{thm:667}. Then, Theorem~\ref{thm:466} (with the colors permuted) implies $|X|\le8$. Since the three blue blocks have total order at most nine, we get $|X|=8$. Theorem~\ref{thm:455}, with the colors permuted, now gives a blue-red or blue-green $K_5$ in $X$. By symmetry, we may assume that it is blue-red. A red neighbor of $X$ exists by minimality and must be a singleton. The other two blue blocks have combined order at most six, giving at most $8+1+6=15$ vertices.

Finally, suppose two blocks $X$ and $Y$ are not all blue. Neither contains a red-green triangle, since such a triangle, a red or green edge in the other block, and one vertex from each blue block would give a red-green $K_7$.

Hence $|X|,|Y|\le6$ by $g(3,6,6)\le7$ (Lemma~\ref{thm:3xy}). The two blue blocks have combined order at most six. To reach seventeen vertices, we need $|X|+|Y|\ge11$, so at least one of $X$ and $Y$ has order six. As in Case~2, Lemma~\ref{lem:sb4} makes each blue block have order at most two (see Figure~\ref{fig:k7-exceptions-step}). The total is, then, at most $6+6+2+2=16$.

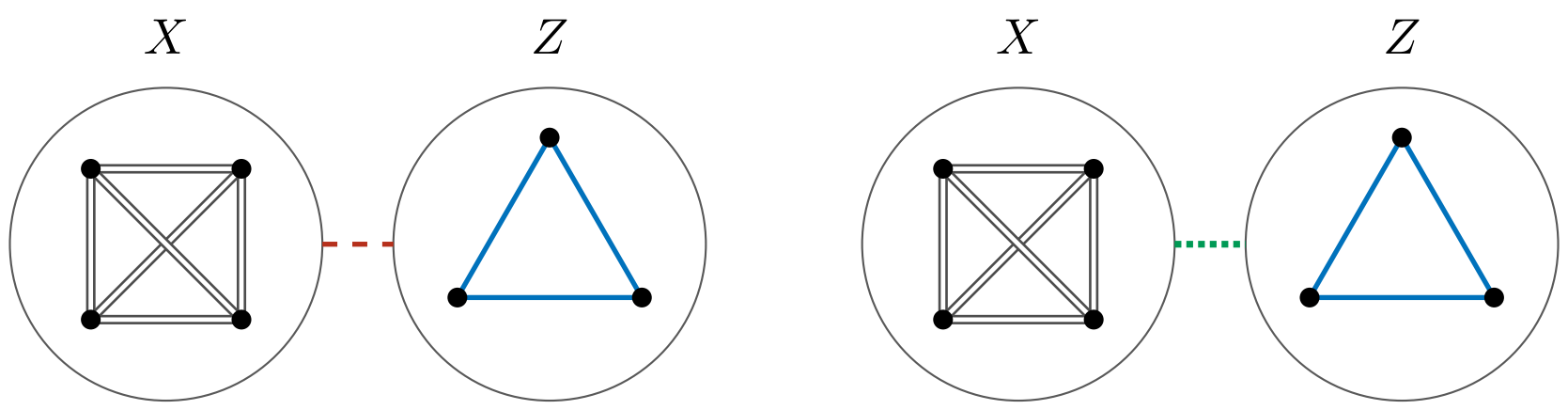
\begin{figure}[htbp]
\centering
\begin{tikzpicture}[every node/.style={font=\small}]
  \foreach \v/\x/\n/\ang/\sty/\lab/\txt in {
    Xr/0/4/45/allowededge/X/{},
    Zr/2.7/3/90/blueedge/Z/{},
    Xg/6/4/45/allowededge/X/{},
    Zg/8.7/3/90/blueedge/Z/{}
  } {
    \node[draw=black!65,circle,minimum size=2.2cm,
          inner sep=0pt,outer sep=0pt] (\v) at (\x,0) {};
    \foreach \i in {1,...,\n} {
      \coordinate (\v-\i) at
        ($(\v.center)+({\ang-360*(\i-1)/\n}:.75)$);
    }
    \foreach \i in {1,...,\n} {
      \foreach \j in {1,...,\n} {
        \ifnum\i<\j\relax
          \draw[\sty] (\v-\i)--(\v-\j);
        \fi
      }
    }
    \foreach \i in {1,...,\n} {
      \node[vertex] at (\v-\i) {};
    }
    \node[above=3pt] at (\v.north) {$\lab$};
    \node[below=3pt,font=\scriptsize,align=center,
          text width=2.2cm] at (\v.south) {\txt};
  }

  \draw[rededge] (Xr)--(Zr);
  \draw[greenedge] (Xg)--(Zg);
\end{tikzpicture}
\caption{$X$ is a block of order six with no red-green triangle. It contains a blue-red $K_4$ (first panel) and a blue-green $K_4$ (second panel). $Z$ is an all-blue block; if $|Z| \geq 3$, the configurations shown are forbidden.}
\label{fig:k7-exceptions-step}
\end{figure}

\case{Case 4: $r=2$}
By symmetry, we may assume that the two blocks $A$ and $B$ are joined in blue. Let $|A|\ge|B|$, and write $a=|A|$ and $b=|B|$. Each block avoids a blue-red and a blue-green $K_6$, since either extends by one outside vertex, and also avoids a red-green $K_7$. Theorem~\ref{thm:667} therefore gives $a\le12$. So, we have $9 \leq a \leq 12$ and $b\ge5$.

If $a\ge11$, then $A$ has a blue-red or blue-green $K_5$ by Theorem~\ref{thm:557}; by symmetry, we may assume that it is blue-red. Thus $B$ is all green, since a blue or red edge there would extend this to a blue-red $K_7$. Since $A$ cannot be all red, it has a blue or green edge. That edge and five vertices of $B$ give a blue-green $K_7$.

The remaining orders are $9+8$ and $10+7$. By $g(4,4,7)=8$ from Theorem~\ref{thm:44n}, the block $A$ contains a blue-red or blue-green $K_4$; by symmetry, we may assume that it is blue-red. Therefore $B$ contains no blue-red triangle, as shown in the first panel of Figure~\ref{fig:k7-final-two-blocks}. Since $b\ge7$ and $B$ has no two-colored $K_7$, it uses all three colors. Then, Corollary~\ref{cor:both4} implies that it has a blue-green $K_4$. Hence, $A$ contains no blue-green triangle (second panel of Figure~\ref{fig:k7-final-two-blocks}).

Finally, $g(3,6,7)\le8$ from Lemma~\ref{thm:3xy}, with the colors permuted, gives a blue-red $K_6$ or a red-green $K_7$ in $A$. The first extends by one vertex of $B$ (third panel of Figure~\ref{fig:k7-final-two-blocks}), and the second is already forbidden. This proves the upper bound.

\begin{figure}[htbp]
\centering
\begin{tikzpicture}[every node/.style={font=\small}]
  \foreach \v/\x/\n/\ang/\sty/\lab/\txt in {
    Aone/0/4/45/allowededge/A/{},
    Bone/2.55/3/90/allowededge/B/{},
    Atwo/5.45/3/90/allowededge/A/{},
    Btwo/8.0/4/45/allowededge/B/{},
    Athree/10.9/6/90/allowededge/A/{}
  } {
    \node[draw=black!65,circle,minimum size=2.2cm,
          inner sep=0pt,outer sep=0pt] (\v) at (\x,0) {};
    \foreach \i in {1,...,\n} {
      \coordinate (\v-\i) at
        ($(\v.center)+({\ang-360*(\i-1)/\n}:.75)$);
    }
    \foreach \i in {1,...,\n} {
      \foreach \j in {1,...,\n} {
        \ifnum\i<\j\relax
          \draw[\sty] (\v-\i)--(\v-\j);
        \fi
      }
    }
    \foreach \i in {1,...,\n} {
      \node[vertex] at (\v-\i) {};
    }
    \node[above=3pt] at (\v.north) {$\lab$};
    \node[below=3pt,font=\scriptsize,align=center,
          text width=2.2cm] at (\v.south) {\txt};
  }

  \node[draw=black!65,circle,minimum size=2.2cm,
        inner sep=0pt,outer sep=0pt] (Bthree) at (13.45,0) {};
  \node[vertex] at (Bthree.center) {};
  \node[above=3pt] at (Bthree.north) {$B$};

  \draw[blueedge] (Aone)--(Bone);
  \draw[blueedge] (Atwo)--(Btwo);
  \draw[blueedge] (Athree)--(Bthree);

  \draw[proofarrow] (3.72,0)--(4.28,0);
  \draw[proofarrow] (9.17,0)--(9.73,0);
\end{tikzpicture}
\caption{The final argument in the proof of Theorem~\ref{thm:k7}, as a chain of contradictions. The first panel has a blue-red $K_4$ in $A$ and a blue-red $K_3$ in $B$. The second panel has a blue-green $K_3$ in $A$ and a blue-green $K_4$ in $B$. The third panel has a blue-red $K_6$ in $A$ and a vertex of $B$.}
\label{fig:k7-final-two-blocks}
\end{figure}
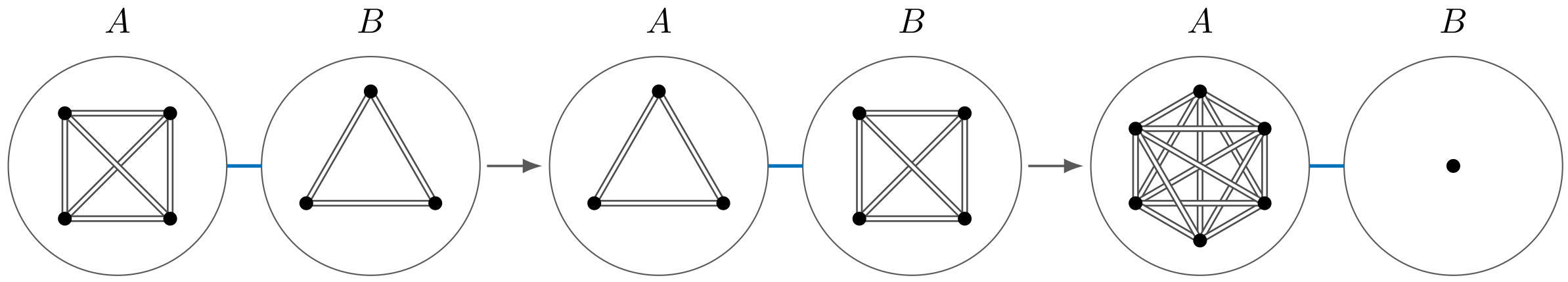

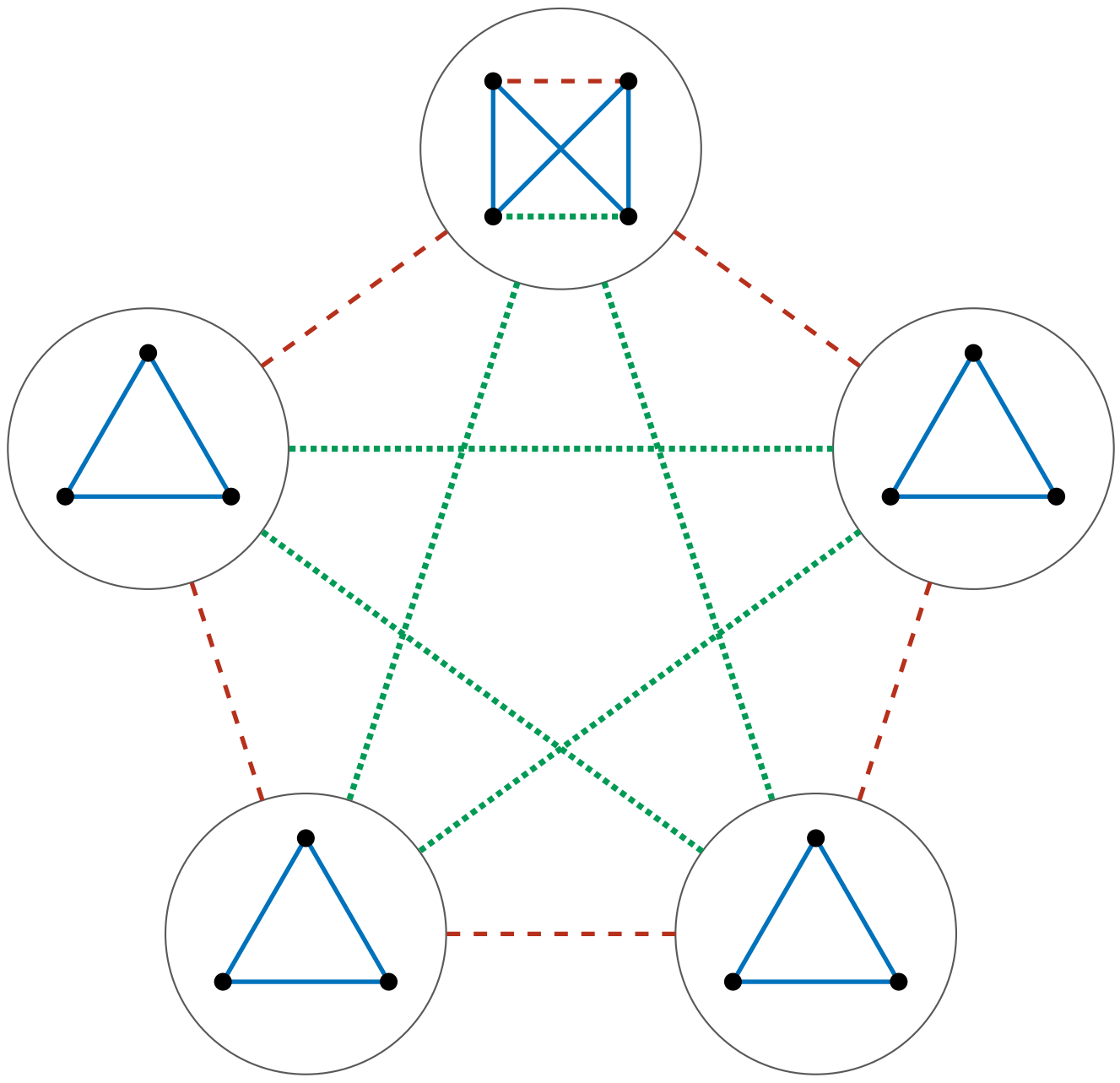
\begin{figure}[htbp]
\centering
\begin{tikzpicture}[every node/.style={font=\small}]
  \foreach \i in {0,...,4} {
    \node[draw=black!65,circle,minimum size=2.2cm,
          inner sep=0pt,outer sep=0pt]
      (B\i) at ({90-72*\i}:3.4) {};
  }
  \foreach \i/\j in {0/1,1/2,2/3,3/4,4/0} {
    \draw[rededge] (B\i)--(B\j);
  }
  \foreach \i/\j in {0/2,0/3,1/3,1/4,2/4} {
    \draw[greenedge] (B\i)--(B\j);
  }

  \foreach \j in {0,...,3} {
    \coordinate (d\j) at
      ($(B0.center)+({45-90*\j}:.75)$);
  }
  \draw[rededge] (d0)--(d3);
  \draw[greenedge] (d1)--(d2);
  \foreach \i/\j in {0/1,0/2,3/1,3/2} {
    \draw[blueedge] (d\i)--(d\j);
  }
  \foreach \j in {0,...,3} {
    \node[vertex] at (d\j) {};
  }
  \foreach \i in {1,...,4} {
    \foreach \j in {0,1,2} {
      \coordinate (b\i-\j) at
        ($(B\i.center)+({90-120*\j}:.75)$);
    }
    \draw[blueedge] (b\i-0)--(b\i-1)--(b\i-2)--cycle;
    \foreach \j in {0,1,2} {
      \node[vertex] at (b\i-\j) {};
    }
  }
\end{tikzpicture}
\caption{A Gallai three-coloring of $K_{16}$ showing $\gr(K_7) \geq 17$. The base graph is $C$; four blocks are blue triangles, and the fifth is a $K_4$ consisting of a blue four-cycle with one red and one green diagonal.}
\label{fig:k16}
\end{figure}

For equality, consider Figure~\ref{fig:k16}. The base graph is the red-green $K_5$ $C$ from Theorem~\ref{thm:456}. Four of the blocks are blue triangles, and the fifth is a $K_4$ whose blue edges form a four-cycle, with one red and one green diagonal.
\end{proof}

\section{AI Disclosure}\label{sec:ai}
The mathematical results in this paper were conceptualized and proved entirely by the authors and did not involve the use of generative AI. Generative AI (GPT 6-Astra, GPT 6.1-Sol, and Claude Opus 5.5) was used to help with minor literature review and proofreading, including spelling and grammar checks and editing the paper. The TikZ code for the figures was also written with the help of AI.

We would like to point out that Lemma~\ref{lem:smallblocks} was suggested as a separate lemma by AI. Before the suggestion, it was implicitly being proved wherever a statement from it was required.

The authors take full responsibility for the work presented in this paper. All errors are attributable to us, not any AI model.

\section*{Acknowledgements}
Gabriel Currier provided extensive assistance, and J\'{o}zsef Solymosi offered guidance and encouragement.

M.A.\ acknowledges that this research was funded in part by the Austrian Science Fund (FWF) \href{https://doi.org/10.55776/AST1535024}{10.55776/AST1535024}. For open access purposes, the authors have applied a CC BY public copyright license to any author accepted manuscript version arising from this submission.

\bibliographystyle{amsplain}
\bibliography{refs}

\end{document}